# DECISION-MAKING: QUALITATIVE INFORMATION

V.Zhukovin , N. Chkhikvadze, Z. Alimbarashvili

**Abstract**
From this set of procedures for given clause we shall choose only interrogation of experts on pairs decisions. It is widely widespread method. It makes the whole chapter in the theory of the decision-making, well investigated with the formal point of view. In the modern theory of decision-making at gathering the initial information for mathematical model it practically does not have alternative.

Incomplete information.

Models of decision-making can be dismembered and presented conditionally in a following kind:

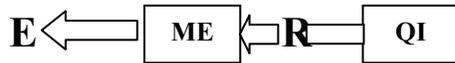

Where **ME** - a mathematical environment of model,
**QI** - Qualitative (psychological) information which is necessary for functioning model and which receive from people. (**PMD** - the person making a decision, experts, the user).
**R**- The means, allowing to connect in a single whole a qualitative and quantitative component of model.
**E** - an output (result).
Possible procedures of reception **QI** the following:
1) **Interrogation**. To the person ask a question which it should answer (including I "do not know").
2) **Physical experiment**. To the person offer some problem situation and fix its behavior, including speech, by means of devices. As an example works V.Pushkina (Moscow) who as the device used game "15", and also R. Nebieridze work (IK Georgia) which fixed movement of eyes at decision problems.
3) **Heuristic experiment** - is association of the first (1) and the second (2) procedures. An example – modeling -fighting of operation for Academy M .Frunze. Certainly, this splitting is very conditional, for in first two procedures demand from the person to prove the answer or behavior, and it already an element of heuristic experiment. Nevertheless it is useful, for allows to define more precisely, we "extort" what information from the person.
Use qualitative information received from the certain contingent of people (experts), in problems of decision-making it is necessary at all stages:
- at formation of set of competitive decisions X,
- at formation of criteria and definition of their scales, or formation Vectors of the Attitude of Preference (VAP) – $R = \{R_j\}_{j=1}^{m}$;

- At an establishment of importance of criteria or $R_j$;
- at a choice of the admissible (comprehensible) conciliatory proposal or ordering X;
- at the account of the additional information.

From this set of procedures for given clause we shall choose only interrogation of experts on pairs decisions. It is widely widespread method. It makes the whole chapter in the theory of the decision-making, well investigated with the formal point of view. In the modern theory of decision-making at gathering the initial information for mathematical model it practically does not have alternative. But it is necessary to use it in any real problem as its weaknesses (by the way it at once are shown concerns " interrogation of experts " as a whole, and not just at pair comparisons).

If the set of concrete decisions contains n alternatives. Number of the ordered pairs equally n (n-1), and disorder -n (n-1) / 2. For n =10 it is 90 and 45. Meanwhile, in practical problems n can be and more than 10. To work with such volume of the information, if still to consider multi-criterial, to any expert it is inconvenient. Therefore, working with the quantitative information, always go by reduction (reduction) of volume of interrogation experts. The example to it is a hierarchical method of the analysis of T.Saaty where it is possible to interrogate the expert on half of general number of pairs. Having theorem Luse $P_{ij} = 1/P_{ji}$ for scales of attitudes, it is possible to work in a scale of differences. All the results presented here , will be in the same scale. Here we shall notice , that all results have analogue in a scale of attitudes.

Let's begin with some definitions:

$Z_{ij} = -Z_{ji}$ - slanting symmetrical function from two alternatives, the analogue of axiom Luse, also defines a degree of the superiority of alternative i above alternative j. As i and j accept values in an interval [1,10], $\|Z_{ij}\|_n$ that is a square matrix, the size nxn. It is analogue of a matrix of T.Saaty in a scale of attitudes.

Now we shall formulate two requirements with which should satisfy probability of interrogation of experts, i.e. the answer.

1) the expert should work with a degree of the superiority $Z_{ij}$. It means, that the expert considers alternative i better (excellent ) alternatives j on any number. And it automatically entails, that the alternative j is worse than alternative i.

2) As a result of full interrogation (at full speed) should receive the linear order, i.e. answers of experts should be coordinated with a condition of transitivity $Z_{ik} + Z_{kj} = Z_{ij}$; and it is analogue of BTL-system. To use this or any other condition of transitivity at interrogation of experts on pair alternatives it is useful.

Let's note, that using only these two conditions, it is possible to fill all matrix, having only one line (or one column). This to term (column) we shall name active line. Thus in a matrix the cross is formed, which will allow the user to fill completely all other cells of a matrix is formed. We shall name it q a cross. (fig.1).

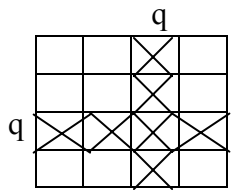

**fig.1**

As $q = 1 \div n$, such crosses will be n.

We shall take an any matrix (i, j). For it $Z_{ij} = Z_{iq} + Z_{qj}$ Both composed are in a cross (i, j).

**The statement 1.** All n crosses are equivalent, as at filling give the same matrix.

**The proof.** At the filled matrix we shall take another p-a cross. It means, that now the p-line is active. We shall take the same cell (i, j). For it we shall write $Z_{ij} = Z_{ip} + Z_{pj}$ But at filling a matrix on the basis of q-a cross had $Z_{ip} = Z_{iq} + Z_{qp}$ and $Z_{pj} = Z_{pq} + Z_{qj}$ Having substituted these formulas in initial for p - a cross, we shall receive the formula for q a-cross, i.e. the statement 1 is proved.

**Consequence 1.** For full filling a matrix $\|Z_{ij}\|$ it is enough to ask the expert to fill one line or one column which becomes thereof active. Thus interrogation meets the requirements formulated above. It is necessary to compare all n pairs alternatives. Not compared pairs are absent. The result of such interrogation - a minus the binary attitude of preference $R_L$, is the linear order $(x_k, x_l) \in R_L$, if $x_k$ is not worse $x_l$. If $(x_k, x_l) \in R_L$ and simultaneously $(x_l, x_k) \in R_L$, then $x_k$ and $x_l$ are equivalent. We Shall clean from $R_L$ set of equivalent pairs $R_L^l$, then we shall deal with coherent strict order. It can be made, as equivalent pairs type $(x_l, x_k) \notin R_L$ bring the zero contribution to comparison of alternatives and a choice of the best. Let's designate the received strict order through $R_L^s$. It is certain on set of set of pairs $E = X \times X$. Let this subset $P \subseteq E$. A condition of transitivity for $R_L^s$ the following:

If $(x_k, x_l) \in R_L^s$ and $(x_l, x_s) \notin R_L^s$ also that $(x_k, x_s) \in R_L^s$ for all three $x_k, x_l, x_s$, form pairs in P. This order can be described by means of a degree of the superiority $Z_{ij}$:
$R_L^s = \{(x_i, x_j) \in E | Z_{ij} > 0, x_i \succ x_j\}$ and $Z_{ij} > 0$ when $x_i \succ x_j$, where $\succ$ means the positive answer of the expert.

We shall enter following designations:

$Z_{ij} > 0$ Through +,

$Z_{ij} < 0$ through -,

$Z_{ij} = 0$ through 0.

Then, using these designations, $R_L^s$ it is possible to present in the form of a matrix $\|Q_{ij}\|_n$, It is a special case of a matrix $\|Z_{ij}\|_n$. For this matrix the condition of transitivity will enter the name as follows: $Z_{ij} = Z_{iq} + Z_{qj}$

1) if $Z_{iq} > 0$ and $Z_{qj} > 0$, that $Z_{ij} > 0$, +;
2) if $Z_{iq} > 0$ and $Z_{qj} = 0$, that $Z_{ij} > 0$, +;
   if $Z_{iq} > 0$ and $Z_{qj} > 0$, that $Z_{ij} > 0$, +;
3) if $Z_{iq} = 0$ and $Z_{qj} = 0$, that $Z_{ij} = 0$, 0;            (2)
4) if $Z_{iq} > 0$ and $Z_{qj} > 0$, that $Z_{ij} > 0$, -;
5) if $Z_{iq} > 0$ and $Z_{qj} > 0$, that $Z_{ij} > 0$, -;
   if $Z_{iq} > 0$ and $Z_{qj} > 0$, that $Z_{ij} > 0$, -;

6) if $Z_{iq} > 0$ and $Z_{qj} > 0$, that $Z_{ij} > 0$, -;

As interrogation spent with performance above the resulted two requirements that a matrix $\|Q_{ij}\|_n$ satisfies with it. We shall remind, that is spoken about full interrogation. If we shall take any other cross in this matrix and we use the formula 2 for the formula 1 we shall receive the same matrix $\|Q_{ij}\|_n$.

**The statement 2.** Having lead private interrogation for one line or a column of a matrix $\|Q_{ij}\|_n$ and using the formula 2 for filling other cells of a matrix we shall receive the information conterminous with the information, received at full interrogation. It means, that the attitude of preference $Q_l^s$ formulated on private interrogations, coincides with full interrogation $R_l^s$.

We use the nonconventional approach to revealing the initial information for model (procedure) of decision-making. In overwhelming number of works the scale of attitudes is used, we use a scale of differences as on experience were convinced, that it is more expressive (is effective). Except for that we actively use the requirement 2 about transitivity which often leads not connecting with a reality of interrogations. Thus, procedure of reception of the full information adequate to full interrogation if it was spent, on the basis of private interrogation looks as follows.

1. We form set of competitive alternatives X. It contains n alternatives $x_i \in X, i = 1 \div n$.
2. We number alternatives from 1 up to $n : \{x_1, x_2, \ldots, x_n\}$.
3. We choose one alternative, we shall tell $x_q \in X$. It means, that active line is q the line.
4. We form n pairs alternatives $(x_q, x_i) \in E$, $i = 1 \div n$,

Considering, that $(x_q, x_q)$ gives 0. We offer the expert for an estimation. Its answers it is entered in q - a line, using signs +, 0, -. If in q - a line is not present - it is the best alternative. Otherwise we shall receive the linear order $R_l(q)$.

5. Using corresponding formulas we shall fill all matrix. Thus we shall receive the linear order - $R_l(P)$ - full. Certainly, both at filling a matrix, and at full interrogation of the expert if it was spent, we carry out both above resulted requirements. Therefore $R_l(P) = R_l$, where $R_l$ - real interrogation. $R_l(P)$ - the linear order. Therefore it Pareto- set consists of one alternative and it the best.

6. Therefore in a matrix there should be one line, shall tell s-concealling, which does not contain any -. In this case $x_s$ - the best alternative. The end of procedure. One of questions at issue is equality $R_l(P) = R_l$.

Let in $R_l$ the pair alternatives $(k,l) \in E$ is estimated to experts positively, i. e., $Z_{kl} > 0$, +. This estimation will keep for all n crosses of full interrogation. Besides, if to lead private interrogation for any cross and to fill a matrix and in it this condition will be carried out. We shall write: $Z_{kl} = Z_{kq} + Z_{ql} > 0$. If the expert will overestimate this pair, i.e. does $Z_{kl} > 0$ -, we shall receive, generally speaking, other linear order. Thus, certainly $R_l(P) \neq R_l^q$. But the matter is that, that $R_l^q$ satisfied to the requirement 2, it is necessary to make changes in initial q-a line. In it either $Z_{iq}$, or $Z_{qj}$, or both together become negative. These at filling a matrix gives the linear order $R_l^{'}(P) = R_l^{'}$, but in any way from this changed line any more we shall not receive the linear order $R_l(P) = R_l$. By the way, and any other order.

**The statement 3.** Private interrogation of the expert on q to a-line, gives the sufficient information for reception of the linear order $R_l(P)$ equivalent real to the linear order $R_l$, received at full interrogation if it was spent. This linear order unique.